\documentclass[review]{elsarticle}

\usepackage{lineno,hyperref}

\usepackage{hyperref}
\usepackage[latin1,utf8x]{inputenc}
\usepackage[english]{babel}
\usepackage{amssymb,amsmath,amsfonts,amsthm}
\usepackage{mathtools}

\usepackage{enumitem}

\theoremstyle{plain}
\newtheorem{teo}{Theorem}[section]

\newtheorem{cor}{Corollary}[section]
\newtheorem{lem}{Lemma}[section]
\newtheorem{proposizione}{Proposition}[section]
\newtheorem*{teo*}{Theorem}
\theoremstyle{definition}
\newtheorem{definizione}{Definition}[section]

\newtheorem{oss}{Remark}[section]
\newtheorem{prop}{Proposition}[section]

\newtheorem*{claim}{Claim}
\newtheorem*{cdimo}{Proof of claim}
\newtheorem{ex}{Example}[section]
\newtheorem{ese}{Exercise}[section]
\newtheorem{conj}{Conjecture}[section]

\newcommand{\bpro}{\begin{proposizione}}		
\newcommand{\epro}{\end{proposizione}}		
\newcommand{\bdef}{\begin{definizione}}		
\newcommand{\bteo}{\begin{teo}}			
\newcommand{\eteo}{\end{teo}}			
\newcommand{\bclaim}{\begin{claim}}			
\newcommand{\eclaim}{\end{claim}}			
\newcommand{\bconj}{\begin{conj}}			
\newcommand{\econj}{\end{conj}}			
\newcommand{\boss}{\begin{oss}}			
\newcommand{\eoss}{\end{oss}}			
\newcommand{\bprop}{\begin{prop}}                        
\newcommand{\eprop}{\end{prop}}                           
\newcommand{\bcor}{\begin{cor}}			
\newcommand{\ecor}{\end{cor}}			
\newcommand{\bdimo}{\begin{proof}}			
\newcommand{\edimo}{\end{proof}}			
\newcommand{\bcdimo}{\begin{cdimo}}			
\newcommand{\ecdimo}{\end{cdimo}}			
\newcommand{\bequ}{\begin{equation}}		
\newcommand{\eequ}{\end{equation}}			
\newcommand{\blem}{\begin{lem}}			
\newcommand{\elem}{\end{lem}}			
\newcommand{\bex}{\begin{ex}}			
\newcommand{\eex}{\end{ex}}				
\newcommand{\bese}{\begin{ese}}			
\newcommand{\eese}{\end{ese}}			
\newcommand{\Z}{\mathbb{Z}}				
\newcommand{\z}{\textbf{\z}}				

\newcommand{\res}{\hat{\phantom{\ }}}

\newcommand{\blom}{\begin{linenomath}}
\newcommand{\elom}{\end{linenomath}}

\begin{document}

\begin{frontmatter}

\title{A generalization of sumsets modulo a prime
}

\author[F. Monopoli]{Francesco Monopoli}       
\address{Dipartimento di Matematica, Universit\`a degli Studi di Milano\\
     via Saldini 50, Milano\\
     I-20133 Italy}
\ead{francesco.monopoli@unimi.it}

     \begin{abstract}
Let $A$ be a set in an abelian group $G$. For integers $h,r \geq 1$ the generalized $h$-fold sumset, denoted by $h^{(r)}A$, is the set of sums of $h$ elements of $A$, where each element appears in the sum at most $r$ times. If $G=\Z$ lower bounds for $|h^{(r)}A|$ are known, as well as the structure of the sets of integers for which $|h^{(r)}A|$ is minimal.
In this paper we generalize this result by giving a lower bound for $|h^{(r)}A|$ when $G=\Z/p\Z$ for a prime $p$, and show new proofs for the direct and inverse problems in $\Z$.

     \end{abstract}


\end{frontmatter}


\begin{section}{Introduction}

Let $A=\{a_1, \dots,  a_k\}$ be a set of $k$ elements in an abelian group $G$. 

Given integers $h, r \geq 1$ define
\blom$$h^{(r)} A = \left\{ \sum_{i=1}^k r_i a_i : 0 \leq r_i \leq r \mbox{ for } i=1, \dots, k \mbox{ and } \sum_{i=1}^k r_i = h\right\}.$$
\elom
Note that the usual sumsets
\blom$$hA = \{ a_{j_1} + \dots + a_{j_h}: a_{j_i} \in A \  \forall i=1, \dots, h\}$$\elom
and the restricted sumsets 
\blom$$h \res A = \{ a_{j_1} + \dots + a_{j_h}: a_{j_i} \in A \ \forall i=1, \dots, h, a_{j_x} \neq a_{j_y} \mbox{ for $x \neq y$}\}$$\elom
can be recovered from this notation, since $hA=h^{(h)} A$ and $ h\res A = h^{(1)} A.$

When $G=\Z$ lower bounds for the cardinality of sumsets and restricted sumsets are well-known.

In this setting, the problem of giving lower bounds for the cardinality of $h^{(r)}A$ for nontrivial values of $h, r$ and $k$ has been studied in \cite{mp}, where the authors proved the following theorem holding for subsets of the integers.

\bteo \label{mpdir}
 Let $h, r$ be non-negative integers, $h=mr + \epsilon, 0 \leq \epsilon \leq r-1$. Let $A$ be a nonempty finite set of integers with $|A|=k$ such that $1 \leq h \leq rk$. Then
\blom\begin{equation}\label{eqfond}|h^{(r)} A| \geq hk - m^2 r +1 -2m\epsilon - \epsilon.\end{equation}\elom
\eteo
Here the condition $h \leq rk$ is necessary, for otherwise the set $h^{(r)}A$ would be empty. 

The lower bound in Theorem \ref{mpdir} is the best one possible, as shown by any arithmetic progression.

A generalization of theorem \ref{mpdir} can be found in \cite{yachen}, where the authors proved lower bounds for generalized sumsets where the $j$th element of $A$ can be repeated up to $r_j$ times, with the $r_j$'s not necessarily all equal to $r$.

In the first section of this paper we will exhibit a new proof of theorem \ref{mpdir}

In the second section we prove the main result of the paper, which states that a similar lower bound also holds when $G=\Z/p\Z$ for a prime $p$.

\bteo \label{mpmodp}
Let $h=mr + \epsilon, 0 \leq \epsilon \leq r-1$. Let $A \subseteq \Z/p \Z$ be a nonempty set with $|A|=k$ such that $1 \leq r \leq h \leq rk$. Then
\blom$$|h^{(r)} A| \geq \min(p, hk - m^2 r +1 -2m\epsilon - \epsilon).$$\elom
\eteo

The authors in \cite{mp} also solved the inverse problem related to $h^{(r)}A$, showing that, up to a few exceptions, any set $A$ satisfying \eqref{eqfond} must be an arithmetic progression:

\bteo\label{mpinv}
Let $k\geq 5$. Let $r$ and $h=mr + \epsilon, 0 \leq \epsilon \leq r-1$ be integers with
 $2 \leq r \leq h \leq rk-2$.
Then any set of $k$ integers $A$ such that
\blom\begin{equation}\label{ugfond}|h^{(r)}A|= hk - m^2 r +1 -2m\epsilon - \epsilon \end{equation}\elom
is a $k$-term arithmetic progression.
\eteo

In the third section we show how we can deduce Theorem \ref{mpinv} from the results in the first two sections and discuss the analogue problem in groups of prime order.
\end{section}
\begin{section}{Direct problem}

To prove Theorem \ref{mpdir}, and later Theorem \ref{mpmodp}, we first deal with the case $r|h$, showing that for a subset $A$ of an abelian group $G$ we have $h^{(r)}A = r(m\res A)$.

\blem\label{lem1}
If $h=mr$, $A \subseteq G$, $|A|=k$ and $rk \geq h \geq 1$
, then $$h^{(r)}A = r(m\res A).$$
\elem

\bdimo
Clearly $r(m\res A) \subseteq h^{(r)} A$, since no element in $A$ can be summed more than $r$ times in order to get an element of $r(m\res A)$.

To prove the converse inclusion, take $x \in h^{(r)}A$ so that, after reordering the elements of $A$ if necessary, $x = \sum_{i=1}^l r_i^{(0)} a_i $ with $ 1 \leq l \leq k, 1 \leq r_i^{(0)} \leq r$ and $\sum_{i=1}^l r_i^{(0)} =h$. 
Let also $r_i^{(0)} = 0$ for $l+1 \leq i \leq k$.

We now describe an algorithm which shows how we can write $x$ as an element in $r(m \res A)$.

If possible, for every $j = 1, \dots, r$ take distinct elements $r_{j_1}^{(j-1)}, \dots, r_{j_m}^{(j-1)}$ which are greater or equal to the remaining $r_{s}^{(j-1)}$ and define
\blom$$x_j = \sum_{i=1}^m a_{j_i}, $$\elom
\blom\begin{equation}\label{eq1} r_{s}^{(j)} = \begin{cases} r_{s}^{(j-1)} -1 & \mbox{ if $s = j_i$ for some $i=1, \dots, m$}\\
		 r_s^{(j-1)} &\mbox{ otherwise}.
\end{cases}
\end{equation}\elom

If we can apply this pocedure for every $j=1,\dots, r$, then we can write $x=x_1 + \dots + x_r$ with $x_i \in m \res A$, thus proving $h^{(r)} A \subseteq r(m \res A)$.

To do this we need to prove that at every step  $j=1, \dots, r$ the following two conditions are satisfied:
\begin{enumerate}
\item\label{con1} $|\{ r_i^{(j-1)} \geq 1 \}_{i}| \geq m$,
\item\label{con2} $\max_{1 \leq i \leq k}(r_i^{(j)}) \leq r-j.$
\end{enumerate}
Since $\sum_{i=1}^k r_i^{(0)} = h = mr$, the first condition holds for $j=1$, and so we can define $r_i^{(1)}$ as in \eqref{eq1}. Clearly $\max_i(r_i^{(1)}) \leq r-1$, for otherwise we could find $m+1$ distinct indexes $s$ such that $r_s^{(0)} = r$, which would imply $\sum_{i=1}^k r_i^{(0)} \geq (m+1)r > h$, a contradiction.

Suppose now that condition \eqref{con1} does not hold for every $j \in [1,r]$, and let $j'$ be the minimal $j$ such that 
\blom	$$|\{ r_i^{(j'-1)} \geq 1 \}| = N < m.$$\elom

By what observed above we must have $2 \leq j' \leq r$.

We have
\blom$$ r_i^{(j'-2)} \begin{cases} 
>1 & \mbox{ for $a$ indexes, $a\leq N < m$} \\
=1 & \mbox{ for $b$ indexes} \\
=0 & \mbox{ for all the remaining $k-a-b$ indexes},
\end{cases}$$\elom
so that $N= a+b -(m-a) = 2a +b-m.$

By the minimality of $j'$ we also have that $a+b \geq m$. 

Next we show that condition \eqref{con2} holds for all $0\leq j'' \leq j'-2 \leq r-2$.

In fact, if this does not happen, take the minimal $j'' \leq j'-2$ which fails to satisfy condition \eqref{con2}, i.e.
\blom$$\max_{1 \leq i \leq k} (r_i^{(j'')}) \geq r-j''+1.$$\elom

By the minimality of $j''$ we must have that $r_i^{(j''-1)} = r- (j''-1)$ for at least $m+1$ values of $i$, because of how the $r_i^{(j)}$ are recursively defined in \eqref{eq1}.

This implies that
\blom$$h-m(j''-1) = \sum_{i=1}^k r_i^{(j''-1)} \geq (m+1)(r-j''+1) = h - m(j''-1) + r-j''+1,$$\elom
a contradiction since $r \geq j''$.

Hence we have that for all $0\leq j'' \leq j'-2$ condition \eqref{con2} is satisfied, which means $\max_i(r_i^{(j'')}) \leq r-j''$.

In particular, since $2a+b = N+m < 2m$ and $a< m$, we get 
\blom\begin{align*}
 h-m(j'-2) &=\sum_{i=1}^k r_i^{(j'-2)} \\
	&\leq a(r-(j'-2)) + b \\
	&<m(r-j')+2m \\
&	= h - m(j'-2),
\end{align*}\elom
a contradiction.

Hence conditions \eqref{con1} and \eqref{con2} are satisfied for all $j=1, \dots, r$.
\edimo

Before proving Theorem \ref{mpdir}, recall the following well-known results on the cardinality of sumsets and restricted sumsets.

\bteo\label{n1} \cite[Theorem 1.3]{nat3}
Let $h\geq 2$.  Let $A$ be a nonempty finite set of integers with $|A|=k$. Then
\blom$$|hA| \geq hk-h+1.$$\elom
\eteo

\bteo\label{n2} \cite[Theorem 1.9]{nat3}
Let $h\geq 2$.  Let $A$ be a nonempty finite set of integers with $|A|=k$. Then
\blom$$|h\res A| = hk-h^2+1$$\elom
\eteo

\bdimo[Proof of Theorem \ref{mpdir}]

Let $A= \{ a_1 < a_2 < \dots < a_k\}$.

The case $\epsilon = 0$ is covered by Lemma \ref{lem1}, since the claim follows from the lower bounds for sumsets and restricted sumsets. 

From now on, assume $\epsilon \geq 1$.

From the condition  $rk \geq h = mr + \epsilon$ we get $k \geq m+1$.

We split the proof in two cases.

{\bfseries Case 1.} $m+\epsilon \leq k$.

In this case it's easy to see the inclusion
\blom\begin{equation*}\label{inc11}B := (r-1)(m \res A) + (m+\epsilon ) \res A \subseteq h^{(r)} A,\end{equation*}\elom
where both the summands are nonempty and $h=(r-1)m + m- \epsilon$.

Then, by Theorems \ref{n1} and \ref{n2} we have
\blom\begin{eqnarray}\label{hreq1}
\nonumber |h^{(r)} A | &=& |B \coprod (h^{(r)}A \setminus B)| \\
		&\geq&hk - m^2r +1 -2m\epsilon - \epsilon^2 + |h^{(r)}A \setminus B|.
\end{eqnarray}\elom

We can now estimate the cardinality of the remaining set observing that
\blom\begin{eqnarray*}\min A = r \sum_{i=1}^m a_i + \epsilon a_{m+1}, \\
\min B = r \sum_{i=1}^m a_i + \sum_{i=m+1}^{m+\epsilon} a_i.\end{eqnarray*}\elom

If we let \blom$$S_{x, y} = r\sum_{i=1}^m a_i + \sum_{i=1}^{x} a_{m+i} + ya_{m+x} + (\epsilon - x - y) a_{m+x+1},$$\elom with $x \in [1, \epsilon -1], y \in [0, \epsilon -x]$, we have $S_{x,y} \in h^{(r)}A$, and
\blom$$
\begin{matrix*}[l]
S_{1, \epsilon -1} &< S_{1, \epsilon -2 } &< S_{1, \epsilon -3}  &< \dots &< S_{1, 0} \\
			&<S_{2, \epsilon - 3} & < S_{2, \epsilon - 4} &< \dots &< S_{2, 0} \\
			&&			 \dots&                &                        & \\
			&&&<S_{\epsilon - 2, 1} & < S_{\epsilon -2, 0} \\
&&&& <S_{\epsilon-1, 0}.
\end{matrix*}
$$\elom

Moreover, all these elements, except for $S_{\epsilon, 0}$ are in $[\min A, \min B -1]$, which gives 
\blom$$ |(h^{(r)}A \setminus B ) \cap [\min A, \min B -1]| \geq \sum_{i=1}^{\epsilon -1} i = \frac{\epsilon^2 - \epsilon}{2}.$$\elom

A symmetric argument gives
\blom$$ |(h^{(r)}A \setminus B ) \cap [\max B +1, \max A]| \geq \sum_{i=1}^{\epsilon -1} i = \frac{\epsilon^2 - \epsilon}{2}.$$\elom

This, combined with equation \eqref{hreq1}, gives the desired lower bound for $|h^{(r)}A|$.

{\bfseries Case 2}: $m+\epsilon> k$.

As already observed in \cite{mp}, we have $|h^{(r)} A| = |(rk - h)^{(r)} A|$.

Then, if $r-1 \leq m+\epsilon$,
\blom$$|h^{(r)}A| = |(r(k-m-1) + (r-\epsilon))^{(r)}A|,$$\elom
and hence we can argue as in the first case to obtain the desired lower bound.

Suppose now $r-1 > m+\epsilon>k$.
Then
\blom\begin{equation*}\label{inc22}B=(m+\epsilon) ((m+1)\res A) + (r-1 - m -\epsilon)(m \res A) \subseteq h^{(r)}A\end{equation*}\elom
and again
\blom\begin{eqnarray}\label{hreq2}
\nonumber |h^{(r)} A | &=& |B \coprod (h^{(r)}A \setminus B)| \\
		&\geq&hk - m^2r +1 -2m\epsilon - \epsilon  - (m^2 + m) + |h^{(r)}A \setminus B|.
\end{eqnarray}\elom

Observe that
\blom\begin{align*}
\min B&=(m+\epsilon) \sum_{i=1}^{m+1} a_{i} + (r-1-m-\epsilon) \sum_{i=1}^m a_i \\
&=(r-1)\sum_{i=1}^m a_i + (m+\epsilon) a_{m+1}, \\
\min A & = r\sum_{i=1}^m a_i + \epsilon a_{m+1}.
\end{align*}\elom

If we let $$T_{x,y} = (r-1) \sum_{i=1}^m a_i + \epsilon a_{m+1} + \sum_{i=x, i \neq y} ^m a_i + x a_{m+1},$$ with $x \in [1, m], y \in [x, m]$, we have $T_{x, y} \in h^{(r)} A$, and

\blom$$
\begin{matrix*}[l]
\min A 		&< T_{1, m} &< T_{1, m-1}  &< \dots &< T_{1, 1} \\
			&<T_{2, m} & < T_{2, m-1} &< \dots &< T_{2, 2} \\
			&			& \dots & & \\
			&&&<T_{m-1, m} & < T_{m-1, m-1} \\
&&&& <T_{m, m}. 
\end{matrix*}
$$\elom

All these elements but $T_{m, m}$ belong to $[\min A, \min B -1]$, which implies
\blom$$ |(h^{(r)}A \setminus B ) \cap [\min A, \min B -1]| \geq \sum_{i=1}^{m} i = \frac{m^2 +m}{2}.$$\elom

A symmetric argument gives
\blom$$ |(h^{(r)}A \setminus B ) \cap [\max B+1, \max A]| \geq \sum_{i=1}^{m} i = \frac{m^2 +m}{2},$$\elom

thus leading, combined with \eqref{hreq2}, to the desired lower bound.
\edimo
\end{section}
\begin{section}{Direct problem in groups of prime order}
In order to prove Theorem \ref{mpmodp} we need the analogues of Theorems \ref{n1} and \ref{n2} in $\Z/p\Z$.

\bteo[Cauchy-Davenport]\label{cd}
Let $h \geq 1$. Let $A \subseteq \Z/p\Z$ be a nonempty set of residues modulo a prime $p$ with $|A|=k$. Then
\blom$$|hA| \geq \min(p, hk-h+1).$$\elom
\eteo

\bteo[Erdős-Heilbronn]\label{eh}
Let $h \geq 1$. Let $A \subseteq \Z/p\Z$ be a nonempty set of residues modulo a prime $p$ with $|A|=k$. Then
\blom$$|h\res A| \geq \min(p, hk-h^2+1).$$\elom
\eteo

Theorem \ref{eh} was conjectured by Erdős and Heilbronn and proved in \cite{dasil} by Da Silva and Hamidoune and later, using the polynomial method, in \cite{pol} by Alon, Nathanson and Ruzsa.

\bdimo[Proof of Theorem \ref{mpmodp}]
The proof goes by induction on $\epsilon$.

If $\epsilon = 0$, thanks to Lemma \ref{lem1} and Theorems \ref{cd} and \ref{eh}, we have:
\blom\begin{eqnarray*}|h^{(r)}A| &=& |r(m \res A)| \geq \min(p, r|m \res A| - r +1) \\
&\geq& \min(p, r\min(p, mk - m^2+1)) - r +1) \\
& = & \min(p, hk - rm^2 +1),
\end{eqnarray*}\elom
where the last equality follows since if $p \leq mk - m^2 +1$ then, for $r\geq1$, 
 $p\leq hk -rm^2+1$.

Let now $\epsilon \in [1, r-1]$.

From $rk \geq h = mr+\epsilon$ we get $k \geq m+1$, and so $h-m-1 = m(r-1) + \epsilon -1 = m(r-1) + \epsilon' \leq (m+1)(r-1) \leq k(r-1)$.

We then have the following inclusion
\blom\begin{equation}\label{inc}(m+1)\res A + (h-m-1)^{(r-1)}A \subseteq h^{(r)}A,\end{equation}\elom
where both summands are nonempty because of the inequalities above.

Moreover, $\epsilon' \in [0, r-2], \epsilon' < \epsilon$ and so, by the inductive hypothesis and Theorems \ref{cd} and \ref{eh} we have
\blom\begin{eqnarray}\label{lasteq}
\nonumber|h^{(r)} A| &\geq& |(m+1)\res A + (h-m-1)^{(r-1)}A| \\
\nonumber&\geq& \min(p, |(m+1) \res A| + |(h-m-1)^{(r-1)}A| -1) \\
&=&\min(p, hk - m^2r -2m\epsilon - \epsilon +1).
\end{eqnarray}\elom
\edimo

Since the inclusion \eqref{inc} holds in any group, our proof, with the obvious modifications, still holds in any abelian group in which theorems similar to \ref{cd} and \ref{eh} hold. See \cite{kar2} for an extensive treatment of the subject.

In particular, when adapted to $\Z$, this leads to yet another proof of Theorem \ref{mpdir}.
\end{section}
\begin{section}{Inverse problem}
From our proof of Theorem \ref{mpdir} it's easy to deduce the inverse theorem based on the well-known results for sumsets and restricted sumsets:

\bteo\label{invn1} \cite[Theorem 1.5]{nat3}
Let $h\geq 2$.  Let $A_1, A_2, \dots, A_h$ be $h$ nonempty finite sets of integers. Then
\blom$$|A_1 + \dots + A_h| = |A_1| + \dots + |A_h|-h+1$$\elom
if and only if the sets $A_1, \dots, A_h$ are arithmetic progressions with the same common difference.
\eteo

\bteo\label{invn2} \cite[Theorem 1.10]{nat3}
Let $h\geq 2$.  Let $A$ be a nonempty finite set of integers with $|A|=k \geq 5$, $2 \leq h \leq k-2$. Then
\blom$$|h \res A| = hk-h^2+1$$\elom
if and only if $A$ is a $k$-term arithmetic progression.
\eteo

\bdimo[Proof of Theorem \ref{mpinv}]
First of all observe that the hypothesis on $h, r$ and $k$ imply that $m \leq k-1$.

Consider first the case $r|h$.

If $m=1$, then $h^{(r)}A=r^{(r)}A=rA$, and Theorem \ref{invn1} can be applied to obtain the thesis.

Let $m \geq 2$. Since $\epsilon = 0$, by Lemma \ref{lem1} we have
\blom$$h(k-m) +1 = |h^{(r)}A| = |r(m \res A)| \geq r|m\res A|-r+1 \geq h(k-m)+1.$$ \elom\

Hence all inequalities above are actually equalities.

In particular, by Theorem \ref{invn1}, $m \res A$ must be an arithmetic progression.

If $m=k-1$, then
\blom\begin{equation}\label{uneq}(k-1)\res A = \left\{ \left(\sum_{i=1}^k a_i\right) - a_k < \left(\sum_{i=1}^k a_i\right) - a_{k-1} < \dots \left(\sum_{i=1}^k a_i\right) - a_1\right\},\end{equation}\elom
and clearly this set is an arithmetic progression if and only if $A$ is an arithmetic progression too.

If $2 \leq m \leq k-2$ we can apply Theorem \ref{invn2} to get the thesis.

Let now $h=mr + \epsilon, \epsilon \in [1, r-1]$.

For $m=0$ we have $h^{(r)}A= \epsilon^{(r)}A = \epsilon A$, and Theorem \ref{invn1} is enough to finish the proof.

Recalling that
$(m+1)\res A + (h-m-1)^{(r-1)}A \subseteq h^{(r)}A,$
from the equation \eqref{ugfond} we deduce that 
\blom\begin{equation}
\label{deq1}|(m+1)\res A| = (m+1)k - (m+1)^2 +1 \end{equation}\elom
and
\blom\begin{equation}
\label{deq2}|(h-m-1)^{(r-1)} A| = (h-m-1)k - m^2(r-1) +1 -2m(\epsilon-1) -\epsilon +1.\end{equation}\elom

By Theorem \ref{invn2} we get the desired conclusion from \eqref{deq1} if $2 \leq m+1 \leq k-2$.

Since we already know that $m+1 \leq k$, only the cases $m=k-2$ and $m=k-1$ are left to study.

If $m=k-2$, then $(m+1)\res A = (k-1) \res A$ and, since \eqref{uneq} holds, we get the thesis.

If $m=k-1$, then $(h-m-1)^{(r-1)} A = (h-k)^{(r-1)}$, and
\blom$$|(h-k)^{(r-1)} A|=|[(r-1)k - h+k]^{(r-1)} A| =|(r-\epsilon)^{(r-1)} A| = |(r-\epsilon) A|$$\elom
since $r-\epsilon \in [1, r-1]$.

This, combined with equation \eqref{deq2} and Theorem \ref{invn1}, gives the desired conclusion.

\edimo

As far as the inverse problem modulo a prime is concerned, in \cite{kar1} the inverse theorem of the  Erdős-Heilbronn conjecture is proved.
\bteo\label{inveh}
Let $A$ be a set of residue classes modulo a prime $p$ with $|A|=k \geq 5, p > 2k-3$. Then
\blom$$|2 \res A| = 2k-3$$\elom
if and only if $A$ is a $k$-term arithmetic progression.
\eteo

The proof however works only when adding two copies of $A$ and, to the best of the author's knowledge, an inverse theorem for $h \res A$, $h > 2$, does not exist yet.

Clearly, an inverse theorem for $h^{(r)} A$ would imply such a result. However, the inclusion \eqref{inc} shows that the converse also holds, showing that the two inverse problems are actually equivalent.

\end{section}

\section*{References}
\bibliography{hfoldbib}
\end{document}